\documentclass[11pt]{article}
\usepackage{amsfonts}
\usepackage{amsmath}

\newtheorem{theorem}{Theorem}
\newtheorem{acknowledgement}[theorem]{Acknowledgement}

\newtheorem{corollary}[theorem]{Corollary}

\newtheorem{lemma}[theorem]{Lemma}

\newtheorem{proposition}[theorem]{Proposition}
\newtheorem{remark}[theorem]{Remark}

\newenvironment{proof}[1][Proof]{\noindent\textbf{#1.} }{\ \rule{0.5em}{0.5em}}
\input{tcilatex}

\begin{document}

\title{On the Weyl Representation of Metaplectic Operators}
\author{Maurice A de Gosson \\
Universit\"{a}t Potsdam, Institut f. Mathematik\\
Am Neuen Palais 10, D-14415 Potsdam\\
maurice.degosson@gmail.com}
\maketitle

\begin{abstract}
We study the Weyl representation of metaplectic operators associated to a
symplectic matrix having no non-trivial fixed point, and justify a formula
suggested in earlier work of Mehlig and Wilkinson. We give precise
calculations of the associated Maslov-type indices; these indices intervene
in a crucial way in Gutzwiller's formula of semiclassical mechanics, and are
simply related to an index defined by Conley and Zehnder.
\end{abstract}

Received 26 December 2004, revised 18 February 2005



\section{Introduction}

In a remarkable paper \cite{MW} Mehlig and Wilkinson propose a simple
derivation of Gutzwiller's \cite{Gutzwiller} approximation%
\begin{equation}
\widetilde{\rho }_{\text{Gutz}}(E)=\frac{1}{\pi \hbar }\func{Re}\sum_{po}%
\frac{T_{po}i^{\nu _{po}}}{\sqrt{|\det (S_{po}-I)|}}e^{i\mathcal{A}%
_{po}/\hbar }  \label{Gutz}
\end{equation}%
for the oscillating part of the semiclassical level density for chaotic
systems whose periodic orbits \textquotedblleft \textit{po}%
\textquotedblright\ are all isolated and non-degenerate ($T_{po}$ is the
prime period, $\nu _{po}$ an integer related to the Maslov index, $\mathcal{A%
}_{po}$ the action, and $S_{po}$ the stability matrix). Mehlig and
Wilkinson's derivation heavily relies upon their observation that for any
symplectic matrix $S$ such that $\det (S-I)\neq 0$ one has%
\begin{equation*}
\frac{i^{\nu }}{\sqrt{|\det (S-I)|}}=\limfunc{Tr}[\widehat{R}_{\nu }(S)]
\end{equation*}%
where $\widehat{R}_{\nu }(S)$ is the operator $L^{2}(\mathbb{R}%
_{x}^{n})\longrightarrow L^{2}(\mathbb{R}_{x}^{n})$ defined by%
\begin{equation*}
\widehat{R}_{\nu }(S)\Psi (x)=\left( \frac{1}{2\pi }\right) ^{n}\frac{i^{\nu
(S)}}{\sqrt{|\det (S-I)|}}\int e^{\frac{i}{2}\left\langle
M_{S}z_{0},z_{0}\right\rangle }\widehat{T}(z_{0})\Psi (x)d^{2n}z_{0}
\end{equation*}%
where:

\begin{itemize}
\item $\widehat{T}(z_{0})$ is the Weyl--Heisenberg operator associated to $%
z_{0}=(x_{0},p_{0})\in \mathbb{R}^{2n}$:%
\begin{equation*}
\widehat{T}(z_{0})f(x)=e^{i(\left\langle p_{0},x\right\rangle -\tfrac{1}{2}%
\left\langle p_{0},x_{0}\right\rangle )}f(x-x_{0})
\end{equation*}%
for any function $f$ defined on $\mathbb{R}^{2n}$;

\item $M_{S}$ is a symmetric matrix, associated to $S$ by the formula%
\begin{equation}
M_{S}=\tfrac{1}{2}J(S+I)(S-I)^{-1}  \label{ms}
\end{equation}%
$I$ being the $2n\times 2n$ identity matrix and $J$ the standard symplectic
matrix;

\item The integer $\nu (S)$ corresponds to a choice of $\arg \det (S-I)$.
\end{itemize}

In what follows we will write Mehlig--Wilkinson's formula as a Bochner
integral%
\begin{equation}
\widehat{R}_{\nu }(S)=\left( \frac{1}{2\pi }\right) ^{n}\frac{i^{\nu (S)}}{%
\sqrt{|\det (S-I)|}}\int e^{\frac{i}{2}\left\langle
M_{S}z_{0},z_{0}\right\rangle }\widehat{T}(z_{0})d^{2n}z_{0}\text{.}
\label{sf3}
\end{equation}

The validity of Mehlig and Wilkinson's derivation requires that --as these
authors claim-- $\widehat{R}_{\nu }(S)$ is one of the two metaplectic
operators $\pm \widehat{S}$ associated with the symplectic matrix $S$. To
sustain this claim the authors refer to previous work of one of the authors 
\cite{Wilkinson}; they also claim that for all $S,S^{\prime }$ such that $%
\det (S-I)\neq 0$, $\det (S^{\prime }-I)\neq 0$ and $\det (SS^{\prime
}-I)\neq 0$ their operators satisfy%
\begin{equation}
\widehat{R}_{\nu ^{\prime \prime }}(SS^{\prime })=\pm \widehat{R}_{\nu }(S)%
\widehat{R}_{\nu ^{\prime }}(S)\text{.}  \label{rr}
\end{equation}

The purpose of this paper is to fully justify Mehlig and Wilkinson's
statements. There are actually several options available. The \textit{a
priori} most direct strategy would be to use Howe's theory \cite{Howe} of
operators with Gaussian kernels (the \textquotedblleft oscillator
semigroup\textquotedblright\ theory: see \cite{Folland} for a review);
another approach would be to use Hannabuss' character theory \cite{Hannabuss}
for contact transformations (also see the follow-up \cite{bohan} to this
paper for interesting applications to star-products). These two methods
however both have, drawbacks. Howe's theory needs quite a lot of technical
prerequisites and would lead here to unnecessarily long calculations; in
addition it would not be very helpful for the study of the sign ambiguity in
(\ref{rr}) since this point is not really addressed in Howe's work. On the
other hand, Hannabuss' machinery is quite abstract (it makes a heavy use of
group character theory) and the use of this approach would perhaps have a
tendency to obscure things. For these reasons we prefer a more
straightforward line of attack, using standard Weyl calculus together with
the theory of the metaplectic group as developed in \cite{AIF}. This
approach moreover has, as we will see, the overwhelming advantage of
producing simple formulae relating the integer $\nu $ in (\ref{sf3}) to the
usual Maslov index of the metaplectic group (we emphasize that $\nu $ is 
\textit{not} the Maslov index!). This is important, because Gutzwiller's
theory has been plagued since its very beginning by the question of how to
calculate the indices $\nu _{po}$ appearing in the trace formula (\ref{Gutz}%
), as witnessed by the abundant literature devoted to this delicate topic
(see \cite{MSdG} where we discuss these issues and give a rather exhaustive
list of references).

\subsection*{Notations}

We denote by $\sigma $ the canonical symplectic form on the phase space $%
\mathbb{R}_{z}^{2n}=\mathbb{R}_{x}^{n}\times \mathbb{R}_{p}^{n}$:%
\begin{equation*}
\sigma (z,z^{\prime })=\left\langle p,x^{\prime }\right\rangle -\left\langle
p^{\prime },x\right\rangle \text{ \ if \ }z=(x,p)\text{, }z^{\prime
}=(x^{\prime }p^{\prime })
\end{equation*}%
that is, in matrix form%
\begin{equation*}
\sigma (z,z^{\prime })=\left\langle Jz,z^{\prime }\right\rangle \text{ \ \ ,
\ \ }J=%
\begin{bmatrix}
0 & I \\ 
-I & 0%
\end{bmatrix}%
\text{.}
\end{equation*}%
The real symplectic group $Sp(n)$ consists of all linear automorphisms $S$
of $\mathbb{R}_{z}^{2n}$ such that $\sigma (Sz,Sz^{\prime })=\sigma
(z,z^{\prime })$ for all $z,z^{\prime }$. It is a connected Lie group and $%
\pi _{1}(Sp(n))$ is isomorphic to $(\mathbb{Z},+)$. $\mathcal{S}(\mathbb{R}%
^{n})$ is the Schwartz space of rapidly decreasing functions on $\mathbb{R}%
^{n}$ and its dual $\mathcal{S}^{\prime }(\mathbb{R}^{n})$ the space of
tempered distributions.

We will denote by $\limfunc{Inert}R$ the number of negative eigenvalues of a
real symmetric matrix $R$.

\section{Prerequisites}

In this Section we briefly recall the main definitions and properties of the
metaplectic group and of Weyl calculus we will need in the rest of this
paper.

\subsection{Standard theory of $Mp(n)$: Review}

The material of this first subsection is quite classical; see for instance 
\cite{Folland,AIF,Cocycles} and the references therein.

Every $\widehat{S}\in Mp(n)$ is the product of two \textquotedblleft
quadratic Fourier transforms\textquotedblright , which are operators $%
\widehat{S}_{W,m}$ defined on $\mathcal{S}(\mathbb{R}^{n})$ by%
\begin{equation}
\widehat{S}_{W,m}f(x)=\left( \tfrac{1}{2\pi i}\right) ^{n/2}i^{m}\sqrt{|\det
L|}\int e^{iW(x,x^{\prime })}f(x^{\prime })d^{n}x^{\prime }  \label{swm1}
\end{equation}%
where $W$ is a quadratic form in the variables $x,x^{\prime }$ of the type%
\begin{equation}
W(x,x^{\prime })=\tfrac{1}{2}\langle Px,x\rangle -\langle Lx,x^{\prime
}\rangle +\tfrac{1}{2}\langle Qx^{\prime },x^{\prime }\rangle  \label{wplq}
\end{equation}%
with $P=P^{T}$, $Q=Q^{T}$, $\det L\neq 0$. The integer $m$
(\textquotedblleft Maslov index\textquotedblright ) appearing in (\ref{swm1}%
) corresponds to a choice of $\arg \det L$:%
\begin{equation*}
m\pi \equiv \arg \det L\text{ \ }\func{mod}2\pi
\end{equation*}%
and to every $W$ there thus corresponds two different choices of $m$ modulo $%
4$: if $m$ is one choice, then $m+2$ is the other, reflecting the fact that $%
Mp(n)$ is a two-fold covering of $Sp(n)$. The projection $\pi
:Mp(n)\longrightarrow Sp(n)$ is entirely specified by the datum of each $\pi
(\widehat{S}_{W,m})$, and we have $\pi (\widehat{S}_{W,m})=S_{W}$ where%
\begin{equation}
(x,p)=S_{W}(x^{\prime },p^{\prime })\Longleftrightarrow p=\partial
_{x}W(x,x^{\prime })\text{ \ and }p^{\prime }=-\partial _{x^{\prime
}}W(x,x^{\prime })\text{.}  \label{leray}
\end{equation}%
Rewriting these conditions in terms of $P,L,Q$ we get $\ p=Px-L^{T}x^{\prime
}$ and $p^{\prime }=Lx-Qx^{\prime }$; solving these equations in $x$ and $p$
yields%
\begin{equation*}
x=L^{-1}(p^{\prime }+Qx^{\prime })\text{ \ and \ \ }p=(PL^{-1}Q-L^{T})x^{%
\prime }+PL^{-1}p^{\prime }
\end{equation*}%
hence the projection $S_{W}$ of $\widehat{S}_{W,m}$ is just the free
symplectic matrix 
\begin{equation}
S_{W}=%
\begin{bmatrix}
L^{-1}Q & L^{-1} \\ 
PL^{-1}Q-L^{T} & PL^{-1}%
\end{bmatrix}
\label{plq}
\end{equation}%
generated by the quadratic form $W$. Note that if conversely $S$ is a free
symplectic matrix 
\begin{equation}
S=%
\begin{bmatrix}
A & B \\ 
C & D%
\end{bmatrix}%
\in Sp(n)\text{ \ , \ }\det B\neq 0  \label{free}
\end{equation}%
then $S=S_{W}$ with $P=DB^{-1}$, $L=B^{-1}$, $Q=B^{-1}A$. Observe that the
free symplectic $2n\times 2n$ matrices form a \textit{dense subset} of $%
Sp(n) $; this property will be used in the proof of Proposition \ref{below}.

The inverse $\widehat{S}_{W,m}^{-1}=(\widehat{S}_{W,m})^{\ast }$ of $%
\widehat{S}_{W,m}$ is the operator $S_{W^{\ast },m^{\ast }}$ where $W^{\ast
}(x,x^{\prime })=-W(x^{\prime },x)$ and $m^{\ast }=n-m$, $\func{mod}4$.

\subsection{Heisenberg--Weyl operators}

The operators $\widehat{T}(z_{0})$ satisfy the metaplectic covariance
formula:%
\begin{equation}
\widehat{S}\widehat{T}(z)=\widehat{T}(Sz)\widehat{S}\text{ \ \ }(S=\pi (%
\widehat{S}))  \label{meco}
\end{equation}%
for every $\widehat{S}\in Mp(n)$ and $z$. In fact, the metaplectic operators
are the only unitary operators (up to a factor in $S^{1}$) satisfying (\ref%
{meco}):

\begin{quote}
\emph{For every }$S\in Sp(n)$ \emph{there exists a unitary transformation }$%
\widehat{U}$ in $L^{2}(\mathbb{R}^{n})$ \emph{satisfying (\ref{meco}) and }$%
\widehat{U}$ \emph{is uniquely determined apart from a constant factor of
modulus one.}
\end{quote}

The Heisenberg--Weyl operators moreover satisfy the relations

\begin{equation}
\widehat{T}(z_{0})\widehat{T}(z_{1})=e^{i\sigma (z_{0},z_{1})}\widehat{T}%
(z_{1})\widehat{T}(z_{0})  \label{noco1}
\end{equation}

\begin{equation}
\widehat{T}(z_{0}+z_{1})=e^{-\tfrac{i}{2}\sigma (z_{0},z_{1})}\widehat{T}%
(z_{0})\widehat{T}(z_{1})  \label{noco2}
\end{equation}%
as is easily seen from their definition.

\subsection{Weyl operators}

Let $a^{w}=a^{w}(x,D)$ be the Weyl operator with symbol $a$ (which we always
assume to belong to some suitable class, allowing the integrals to be viewed
as distributions): 
\begin{equation*}
a^{w}f(x)=\left( \tfrac{1}{2\pi }\right) ^{n}\diint e^{i\left\langle
p,x-y\right\rangle }a(\tfrac{1}{2}(x+y),p)f(y)d^{n}yd^{n}p\text{; }
\end{equation*}%
where $f\in \mathcal{S}(\mathbb{R}^{n})$; equivalently%
\begin{equation*}
a^{w}=\left( \tfrac{1}{2\pi }\right) ^{n}\int a_{\sigma }(z_{0})\widehat{T}%
(z_{0})d^{2n}z_{0}
\end{equation*}%
where the \textquotedblleft twisted symbol\textquotedblright\ $a_{\sigma }$
is the symplectic Fourier transform $F_{\sigma }a$ of $a$:%
\begin{equation*}
a_{\sigma }(z)=F_{\sigma }a(z)=\left( \tfrac{1}{2\pi }\right) ^{n}\int
e^{-i\sigma (z,z^{\prime })}a(z^{\prime })d^{2n}z^{\prime }\text{.}
\end{equation*}%
The compose $c^{w}=a^{w}\circ b^{w}$ (when defined) is the Weyl operator
with twisted Weyl symbol 
\begin{equation}
c_{\sigma }=\left( \tfrac{1}{2\pi }\right) ^{n}(a_{\sigma }\ast _{\sigma
}b_{\sigma })  \label{cw}
\end{equation}%
where 
\begin{equation}
a\ast _{\sigma }b(z)=\int e^{\frac{i}{2}\sigma (z,u)}a(z-u)b(u)d^{2n}u
\label{tw}
\end{equation}%
(see for instance Littlejohn \cite{Littlejohn}, Wong \cite{Wong}).

\subsection{Generalized Fresnel Formula}

We will use the following formula, generalizing the usual Fresnel integral
to complex Gaussians. Let $M$ be a real symmetric $m\times m$ matrix. If $M$
is invertible then the Fourier transform of the exponential $\exp
(i\left\langle Mx,x\right\rangle /2)$ is given by the formula%
\begin{equation}
\left( \tfrac{1}{2\pi }\right) ^{m/2}\int e^{-i\left\langle v,u\right\rangle
}e^{\frac{i}{2}\left\langle Mu,u\right\rangle }d^{m}u=|\det M|^{-1/2}e^{%
\frac{i\pi }{4}\limfunc{sgn}M}e^{-\frac{i}{2}\left\langle
M^{-1}v,v\right\rangle }  \label{fres}
\end{equation}%
where $\limfunc{sgn}M$, the \textquotedblleft signature\textquotedblright\
of $M$, is the number of $>0$ eigenvalues of $M$ minus the number of $<0$
eigenvalues.

For a proof of this formula see for instance \cite{Folland}, Appendix A.

\section{Discussion of the Mehlig--Wilkinson Formula}

In this Section we show that the Mehlig--Willkinson operators (\ref{sf3})
indeed are metaplectic operators. We begin by giving two alternative
expressions for these operators.

\subsection{Equivalent formulations}

We begin by remarking that the matrix 
\begin{equation*}
M_{S}=\frac{1}{2}J(S+I)(S-I)^{-1}
\end{equation*}
is symmetric; this immediately follows from the conditions%
\begin{equation*}
S\in Sp(n)\Longleftrightarrow S^{T}JS=J\Longleftrightarrow SJS^{T}=J\text{.}
\end{equation*}%
Notice that for every $M$ with $\det (M-\frac{1}{2}J)\neq 0$ the equation%
\begin{equation*}
M=\tfrac{1}{2}J(S+I)(S-I)^{-1}
\end{equation*}%
can be solved\ in $S$, yielding 
\begin{equation*}
S=(M-\tfrac{1}{2}J)^{-1}(M+\tfrac{1}{2}J);
\end{equation*}%
the relation $S\in Sp(n)$ is then equivalent to $M$ being real and symmetric.

\begin{lemma}
Let $S\in Sp(n)$ be such that $\det (S-I)\neq 0$. The operator 
\begin{equation}
\widehat{R}_{\nu }(S)=\left( \frac{1}{2\pi }\right) ^{n}\frac{i^{\nu }}{%
\sqrt{|\det (S-I)|}}\int e^{\frac{i}{2}\left\langle M_{S}z,z\right\rangle }%
\widehat{T}(z)d^{2n}z  \label{alf0}
\end{equation}%
can be written as 
\begin{equation}
\widehat{R}_{\nu }(S)=\left( \frac{1}{2\pi }\right) ^{n}i^{\nu }\sqrt{|\det
(S-I)|}\int e^{-\frac{i}{2}\sigma (Sz,z)}\widehat{T}((S-I)z)d^{2n}z
\label{alf2}
\end{equation}%
that is, as 
\begin{equation}
\widehat{R}_{\nu }(S)=\left( \frac{1}{2\pi }\right) ^{n}i^{\nu }\sqrt{|\det
(S-I)|}\int \widehat{T}(Sz)\widehat{T}(-z)d^{2n}z\text{.}  \label{alf1}
\end{equation}
\end{lemma}

\begin{proof}
We have 
\begin{equation*}
\tfrac{1}{2}J(S+I)(S-I)^{-1}=\tfrac{1}{2}J+J(S-I)^{-1}
\end{equation*}%
hence, in view of the antisymmetry of $J$,%
\begin{equation*}
\left\langle M_{S}z,z\right\rangle =\left\langle J(S-I)^{-1}z,z\right\rangle
=\sigma ((S-I)^{-1}z,z)
\end{equation*}%
Performing the change of variables $z\longmapsto (S-I)z$ we can rewrite the
integral in the right-hand side of (\ref{alf0}) as%
\begin{eqnarray*}
\int e^{\frac{i}{2}\left\langle M_{S}z,z\right\rangle }\widehat{T}(z)d^{2n}z
&=&\sqrt{|\det (S-I)|}\int e^{\frac{i}{2}\sigma (z,(S-I)z)}\widehat{T}%
((S-I)z)d^{2n}z \\
&=&\sqrt{|\det (S-I)|}\int e^{-\frac{i}{2}\sigma (Sz,z)}\widehat{T}%
((S-I)z)d^{2n}z
\end{eqnarray*}%
hence (\ref{alf2}). Taking into account the relation (\ref{noco2}) we have%
\begin{equation*}
\widehat{T}((S-I)z)=e^{\tfrac{i}{2}\sigma (Sz,z)}\widehat{T}(Sz)\widehat{T}%
(-z)
\end{equation*}%
and formula (\ref{alf1}) follows.
\end{proof}

\begin{remark}
Formulae (\ref{alf0}) and (\ref{alf2})--(\ref{alf1}) suggest that $\nu \pi $
could be a choice of $\pm \arg \det (S-I)$. This is however not the case
(see (\ref{G}) in Proposition \ref{prim} below); formula (\ref{cz1}) will
identify the integer $\nu $ with the Conley--Zehnder index.
\end{remark}

\begin{corollary}
\label{cor2}We have $\widehat{R}_{\nu }(S)=c_{S}\widehat{S}_{W,m}$ with $%
|c_{S}|=1$.
\end{corollary}

\begin{proof}
The operator $\widehat{R}_{\nu }(S)$ satisfies the metaplectic covariance
relation 
\begin{equation*}
\widehat{R}_{\nu }(S)\widehat{T}(z)=\widehat{T}(Sz)\widehat{R}_{\nu }(S)
\end{equation*}%
as immediately follows from the alternative form (\ref{alf2}) of $\widehat{R}%
_{\nu }(S)$. On the other hand, a straightforward calculation using formula (%
\ref{alf1}) shows that $\widehat{R}_{\nu }(S)$ is unitary, hence the claim.
\end{proof}

Let us precise Corollary \ref{cor2} by discussing the choice of the constant 
$c_{S}$.

\subsection{The case $\widehat{S}=\widehat{S}_{W,m}$}

We are going to show that the Mehlig--Wilkinson operators coincide with the
metaplectic operators $\widehat{S}_{W,m}$ when $S=S_{W}$ and we will
thereafter determine the correct choice for $\nu $; we will see that it is
related by a simple formula to the usual Maslov index as defined in \cite%
{AIF}.

Let us first prove the following technical result:

\begin{lemma}
\label{lemma1}Let $S_{W}$ be a free symplectic matrix (\ref{free}). We have 
\begin{equation}
\det (S_{W}-I)=(-1)^{n}\det (B)\det (B^{-1}A+DB^{-1}-B^{-1}-(B^{T})^{-1})
\label{bofor1}
\end{equation}%
that is, when $S$ is written in the form (\ref{plq}):%
\begin{equation}
\det (S_{W}-I)=(-1)^{n}\det (L^{-1})\det (P+Q-L-L^{T})\text{.}
\label{bofor2}
\end{equation}
\end{lemma}

\begin{proof}
Since $B$ is invertible we can write $S-I$ as%
\begin{equation}
\begin{bmatrix}
A-I & B \\ 
C & D-I%
\end{bmatrix}%
=%
\begin{bmatrix}
0 & B \\ 
I & D-I%
\end{bmatrix}%
\begin{bmatrix}
C-(D-I)B^{-1}(A-I) & 0 \\ 
B^{-1}(A-I) & I%
\end{bmatrix}
\label{F}
\end{equation}%
and hence%
\begin{equation*}
\det (S_{W}-I)=\det (-B)\det (C-(D-I)B^{-1}(A-I))\text{.}
\end{equation*}%
Since $S$ is symplectic we have $C-DB^{-1}A=-(B^{T})^{-1}$ (use for instance
the fact that $S^{T}JS=SJS^{T}=J$) and hence%
\begin{equation*}
C-(D-I)B^{-1}(A-I))=B^{-1}A+DB^{-1}-B^{-1}-(B^{T})^{-1}\text{;}
\end{equation*}%
the Lemma follows since $\det (-B)=(-1)^{n}\det B$.
\end{proof}

\begin{remark}
The factorization (\ref{F}) shows in particular that $\ker (S-I)$ is
isomorphic to $\ker (P+Q-L-L^{T})$ (cf. \cite{MPT}, Lemma 2.8, and proof of
Lemma 2.9).
\end{remark}

Let us denote by $W_{xx}$ the Hessian matrix of the function $x\longmapsto
W(x,x)$, that is%
\begin{equation*}
W_{xx}=P+Q-L-L^{T}\text{. }
\end{equation*}%
We have:

\begin{proposition}
\label{prim}Let $S=S_{W}$ be a free symplectic matrix (\ref{free}) and $%
\widehat{R}_{\nu }(S)$ the corresponding Mehlig--Wilkinson operator. We have 
$\widehat{R}_{\nu }(S)=\widehat{S}_{W,m}$ provided that $\nu $ is chosen so
that%
\begin{equation}
\nu \equiv m-\limfunc{Inert}W_{xx}\text{ \ }\func{mod}4  \label{Maslov1}
\end{equation}%
in which case we have 
\begin{equation}
\frac{1}{\pi }\arg \det (S-I)\equiv -\nu +n\text{\ \ }\func{mod}2  \label{G}
\end{equation}
\end{proposition}

\begin{proof}
Recall that we have shown that $\widehat{R}_{\nu }(S)=c_{S}\widehat{S}_{W,m}$
where $c_{S}$ is some complex constant with $|c_{S}|=1$. Let us determine
that constant. Let $\delta \in \mathcal{S}^{\prime }(\mathbb{R}^{n})$ be the
Dirac distribution centered at $x=0$; setting%
\begin{equation*}
C_{W,\nu }=\left( \frac{1}{2\pi }\right) ^{n}\frac{i^{\nu }}{\sqrt{|\det
(S-I)|}}
\end{equation*}%
we have, by definition of $\widehat{R}_{\nu }(S)$, writing $%
z_{0}=(x_{0},p_{0})$ in place of $z=(x,p)$: 
\begin{eqnarray*}
\widehat{R}_{\nu }(S)\delta (x) &=&C_{W,\nu }\int e^{\frac{i}{2}\left\langle
M_{S}z_{0},z_{0}\right\rangle }e^{i(\left\langle p_{0},x\right\rangle -\frac{%
1}{2}\left\langle p_{0},x_{0}\right\rangle )}\delta (x-x_{0})d^{2n}z_{0} \\
&=&C_{W,\nu }\int e^{\frac{i}{2}\left\langle
M_{S}(x,p_{0}),(x,p_{0})\right\rangle }e^{\frac{i}{2}\left\langle
p_{0},x\right\rangle }\delta (x-x_{0})d^{2n}z_{0}
\end{eqnarray*}%
hence, setting $x=0$,%
\begin{equation*}
\widehat{R}_{\nu }(S)\delta (0)=C_{W,\nu }\int e^{\frac{i}{2}\left\langle
M_{S}(0,p_{0}),(0,p_{0})\right\rangle }\delta (-x_{0})d^{2n}z_{0}
\end{equation*}%
that is, since $\int \delta (-x_{0})d^{n}x_{0}=1$,%
\begin{equation}
\widehat{R}_{\nu }(S)\delta (0)=\left( \frac{1}{2\pi }\right) ^{n}\frac{%
i^{\nu }}{\sqrt{|\det (S-I)|}}\int e^{\frac{i}{2}\left\langle
M_{S}(0,p_{0}),(0,p_{0})\right\rangle }d^{n}p_{0}\text{.}  \label{sdo}
\end{equation}%
Let us calculate the scalar product 
\begin{equation*}
\left\langle M_{S}(0,p_{0}),(0,p_{0})\right\rangle =\sigma
((S-I)^{-1}0,p_{0}),(0,p_{0}))\text{.}
\end{equation*}%
The relation $(x,p)=(S-I)^{-1}(0,p_{0})$ is equivalent to $%
S(x,p)=(x,p+p_{0}) $ that is to%
\begin{equation*}
p+p_{0}=\partial _{x}W(x,x)\text{ \ and \ }p=-\partial _{x^{\prime }}W(x,x)%
\text{.}
\end{equation*}%
Using the explicit form (\ref{wplq}) of $W$ together with Lemma \ref{lemma1}
these relations yield%
\begin{equation*}
x=(P+Q-L-L^{T})^{-1}p_{0}\text{ \ and \ }p=(L-Q)(P+Q-L-L^{T})^{-1}p_{0}
\end{equation*}%
and hence%
\begin{equation}
\left\langle M_{S}(0,p_{0}),(0,p_{0})\right\rangle =-\left\langle
W_{xx}^{-1}p_{0},p_{0}\right\rangle \text{.}  \label{bofor3}
\end{equation}%
Applying Fresnel's formula (\ref{fres}) we get%
\begin{equation*}
\left( \tfrac{1}{2\pi }\right) ^{n}\int e^{\frac{i}{2}\left\langle
M_{S}(0,p_{0}),(0,p_{0})\right\rangle }d^{n}p_{0}=\left( \tfrac{1}{2\pi }%
\right) ^{n/2}e^{-\frac{i\pi }{4}\limfunc{sgn}W_{xx}}|\det W_{xx}|^{1/2}%
\text{;}
\end{equation*}%
noting that 
\begin{equation*}
\frac{1}{\sqrt{|\det (S_{W}-I)|}}=|\det L|^{1/2}|\det W_{xx}|^{-1/2}
\end{equation*}%
(formula (\ref{bofor2}) in Lemma \ref{lemma1}) we thus have%
\begin{equation*}
\widehat{R}_{\nu }(S_{W})\delta (0)=\left( \tfrac{1}{2\pi }\right)
^{n/2}i^{\nu }e^{-\frac{i\pi }{4}\limfunc{sgn}W_{xx}}|\det L|^{1/2}\text{.}
\end{equation*}%
Now, by definition (\ref{swm1}) of $\widehat{S}_{W,m}$,%
\begin{equation*}
\widehat{S}_{W,m}\delta (0)=\left( \tfrac{1}{2\pi }\right)
^{n/2}i^{m-n/2}|\det L|^{1/2}
\end{equation*}%
hence $i^{\nu }e^{-\frac{i\pi }{4}\limfunc{sgn}W_{xx}}=i^{m-n/2}$. It
follows that we have%
\begin{equation*}
\nu -\tfrac{1}{2}\limfunc{sgn}W_{xx}\equiv m-\tfrac{1}{2}n\text{ \ }\func{mod%
}4
\end{equation*}%
which is the same thing as (\ref{Maslov1}) since $W_{xx}$ has rank $n$. In
view of (\ref{bofor1}) we have%
\begin{equation*}
\tfrac{1}{\pi }\arg \det (S_{W}-I)=n+m+\arg \det W_{xx}\text{ \ }\func{mod}2;
\end{equation*}%
formula (\ref{G}) follows using (\ref{Maslov1}).
\end{proof}

Let us digress for a while on the integers $m$ and $\limfunc{Inert}W_{xx}$
appearing in formula (\ref{Maslov1}) and discuss them from the point of view
of calculus of variations. It is for this purpose useful to recall that in
Gutzwiller's formula (of which Mehlig and Wilkinson precisely want to give a
new approach using the operators $\widehat{R}_{\nu }(S)$) the symplectic
matrix $S$ is obtained from the monodromy matrix of an isolated Hamiltonian
periodic orbit. Let us go a little bit further. Consider a Hamiltonian flow $%
\phi _{t}$ determined by some time-dependent Hamiltonian $H=H(z,t)$ defined
on $\mathbb{R}_{z}^{2n}\times \mathbb{R}_{t},$ and let $z_{0}=(x_{0},p_{0})$
be such that $\phi _{T}(z_{0})=z_{0}$ for some $T>0$. The Jacobian matrices $%
S(z_{0},t)=D\phi _{t}(z_{0})$ are symplectic and satisfy the
\textquotedblleft variational equation\textquotedblright 
\begin{equation*}
\frac{d}{dt}S(z_{0},t)=JH^{\prime \prime }(z_{0},t)S(z_{0},t)
\end{equation*}%
where $H"(z_{0},t)$ is the Hessian matrix $D^{2}H(\phi _{t}(z_{0}),t)$. When 
$t$ varies from $0$ to $T$ the matrices $S(z_{0},t)$ describe a path in $%
Sp(n)$ originating at the identity and ending at $S(z_{0},T)$ (the
\textquotedblleft monodromy matrix\textquotedblright ). Suppose that $%
S(z_{0},T)$ is a free symplectic matrix $S_{W}$; then there exist $p$ and $%
p^{\prime }$ such that $(x_{0},p)=S_{W}(x_{0},p^{\prime }),$ that is,
expressing $W$ in terms of $P,L,Q$ as in (\ref{plq}) and using (\ref{leray}%
), $p=(P-L^{T})x_{0}$ and $p^{\prime }=(L-Q)x_{0}$, that is%
\begin{equation*}
p-p^{\prime }=(P+Q-L-L^{T})x_{0}\text{.}
\end{equation*}%
It thus appears (see for instance \cite{MPT,MG}) that $\limfunc{Inert}W_{xx}$
is Morse's \cite{Morse} \textit{order of concavity} of the periodic orbit
through $z_{0}$\textit{. }

\begin{remark}
In \cite{MPT} Piccione and his collaborators use the order of concavity,
which they identify with $\limfunc{Inert}W_{xx}$, to investigate the Maslov
and Morse indices for periodic geodesics. The topic is also discussed at
some length in Muratore--Ginnaneschi \cite{MG} in connection with the study
of Gutzwiller's formulae using field-theoretical methods.
\end{remark}

Perhaps even more interesting is the relation between the index $\nu $ and
an index defined by Conley and Zehnder in \cite{CZ}. Let us denote by $%
Sp_{0}(n)$ the set of all $S\in Sp(n)$ such that $\det (S-I)\neq 0$. We have%
\begin{equation*}
Sp_{0}(n)=Sp_{+}(n)\cup Sp_{-}(n)
\end{equation*}%
where $S\in Sp_{\pm }(n)$ $if$ and only if $\pm \det (S-I)>0$. The sets $%
Sp_{\pm }(n)$ are connected and every loop in $Sp_{0}(n)$ is contractible in 
$Sp(n)$. Consider now a path $\tilde{S}:[0,T]\longmapsto Sp_{0}(n)$ going
from the identity to $S=S(T)$ ($S(T)$ may be viewed, if one wants, as the
monodromy matrix of a periodic Hamiltonian orbit corresponding to a \textit{%
time-dependent} Hamiltonian). The index of Conley and Zehnder associates to
the path $\tilde{S}$ an integer $\mu _{CZ}(\tilde{S})$ only depending on the
homotopy class (with fixed endpoints) of that path, and such that%
\begin{equation*}
\func{sign}\det (S-I)=(-1)^{n-\mu _{CZ}(\tilde{S})}
\end{equation*}%
that is, equivalently,%
\begin{equation*}
\frac{1}{\pi }\arg \det (S-I)\equiv n-\mu _{CZ}(\tilde{S})\text{ }\func{mod}2%
\text{.}
\end{equation*}%
It follows from formula (\ref{G}) in Proposition \ref{prim} that we have%
\begin{equation}
\nu \equiv \mu _{CZ}(\tilde{S})\text{ }\func{mod}2  \label{cz1}
\end{equation}%
and from formula (\ref{Maslov1}) in the same proposition that 
\begin{equation}
\mu _{CZ}(\tilde{S})\equiv m+n-\limfunc{Inert}W_{xx}\ \func{mod}2\text{.}
\label{cz2}
\end{equation}

\subsection{The general case}

Recall that we established in Lemma \ref{lemma1} the equality%
\begin{equation}
\det (S_{W}-I)=(-1)^{n}\det L^{-1}\det (P+Q-L-L^{T}).  \label{splq}
\end{equation}%
valid for all free matrices $S_{W}\in Sp(n)$. Also recall that every $%
\widehat{S}\in Mp(n)$ can be written (in infinitely many ways) as a product $%
\widehat{S}=\widehat{S}_{W,m}\widehat{S}_{W^{\prime },m^{\prime }}$. We are
going to show that $\widehat{S}_{W,m}$ and $\widehat{S}_{W^{\prime
},m^{\prime }}$ always can be chosen such that $\det (S_{W}-I)\neq 0$ and $%
\det (S_{W^{\prime }}-I)\neq 0$. For that purpose we need the following
straightforward factorization result, which we nevertheless glorify by
putting it into italics:

\begin{lemma}
Let $W$ be given by (\ref{wplq}); then 
\begin{equation}
\widehat{S}_{W,m}=\widehat{V}_{-P}\widehat{M}_{L,m}\widehat{J}\widehat{V}%
_{-Q}  \label{fac1}
\end{equation}%
where 
\begin{equation*}
\widehat{V}_{-P}f(x)=e^{\frac{i}{2}\left\langle Px,x\right\rangle }f(x)\text{
\ , \ }\widehat{M}_{L,m}f(x)=i^{m}\sqrt{|\det L|}f(Lx)\text{ ,}
\end{equation*}%
and $\widehat{J}$ is the modified Fourier transform given by%
\begin{equation*}
\widehat{J}f(x)=\left( \tfrac{1}{2\pi i}\right) ^{n/2}\int e^{-i\left\langle
x,x^{\prime }\right\rangle }f(x^{\prime })d^{n}x^{\prime }\text{.}
\end{equation*}
\end{lemma}

\begin{proof}
It is obvious using the explicit expression (\ref{wplq}) of the quadratic
form $W$ (see \cite{AIF}).
\end{proof}

Let us now state and prove the first result of this section:

\begin{proposition}
\label{above}(i) Every $\widehat{S}\in Mp(n)$ can be written as a product 
\begin{equation}
\widehat{S}=\widehat{R}_{\nu }(S_{W})\widehat{R}_{\nu ^{\prime
}}(S_{W^{\prime }})\text{.}  \label{rsw}
\end{equation}%
(ii) The Mehlig--Wilkinson operators thus generate $Mp(n)$.
\end{proposition}

\begin{proof}
\textit{(ii)} follows from \textit{(i)} since the $\widehat{S}_{W,m}$
generate $Mp(n)$. To prove \textit{(i)} let us write $\widehat{S}=\widehat{S}%
_{W,m}\widehat{S}_{W^{\prime },m^{\prime }}$ and apply (\ref{fac1}) to each
of the factors; letting $P^{\prime },L^{\prime },Q^{\prime }$ define $%
W^{\prime }$ just as (\ref{wplq}) is defined by $P,L,Q$ we have%
\begin{equation}
\widehat{S}=\widehat{V}_{-P}\widehat{M}_{L,m}\widehat{J}\widehat{V}%
_{-(P^{\prime }+Q)}\widehat{M}_{L^{\prime },m^{\prime }}\widehat{J}\widehat{V%
}_{-Q^{\prime }}\text{.}  \label{sprod}
\end{equation}%
We claim that $\widehat{S}_{W,m}$ and $\widehat{S}_{W^{\prime },m^{\prime }}$
can be chosen in such a way that $\det (S_{W}-I)\neq 0$ and $\det
(S_{W^{\prime }}-I)\neq 0$ that is, 
\begin{equation*}
\det (P+Q-L-L^{T})\neq 0\text{ \ and \ }\det (P^{\prime }+Q^{\prime
}-L^{\prime }-L^{\prime T})\neq 0\text{;}
\end{equation*}%
this will prove the assertion in view of (\ref{splq}). We first remark that
the right hand-side of (\ref{sprod}) obviously does not change if we replace 
$P^{\prime }$ by $P^{\prime }+\lambda $ and $Q$ by $Q-\lambda $ where $%
\lambda \in \mathbb{R}$. Choose now $\lambda $ such that it is not an
eigenvalue of $P+Q-L-L^{T}$ and $-\lambda $ is not an eigenvalue of $%
P^{\prime }+Q^{\prime }-L^{\prime }-L^{\prime T}$; then 
\begin{equation*}
\det (P+Q-\lambda I-L-L^{T})\neq 0\text{ \ and \ }\det (P^{\prime }+\lambda
+Q^{\prime }-L-L^{T})\neq 0
\end{equation*}%
and we have $\widehat{S}=\widehat{S}_{W_{1},m_{1}}\widehat{S}_{W_{1}^{\prime
},m_{1}^{\prime }}$ with%
\begin{eqnarray*}
W_{1}(x,x^{\prime }) &=&\tfrac{1}{2}\langle Px,x\rangle -\langle
Lx,x^{\prime }\rangle +\tfrac{1}{2}\langle (Q-\lambda )x^{\prime },x^{\prime
}\rangle \\
W_{1}^{\prime }(x,x^{\prime }) &=&\tfrac{1}{2}\langle (P^{\prime }+\lambda
)x,x\rangle -\langle L^{\prime }x,x^{\prime }\rangle +\tfrac{1}{2}\langle
Q^{\prime }x^{\prime },x^{\prime }\rangle \text{.}
\end{eqnarray*}
\end{proof}

So far, so good. But we haven't told the whole story yet: there remains to
prove that $\widehat{S}\in Mp(n)$ can be written in the form $\widehat{R}%
_{\nu }(S)$ if $\det (S-I)\neq 0$.

\begin{proposition}
\label{below}Let $\widehat{S}\in Mp(n)$ be such that $\det (S-I)\neq 0$. If $%
\widehat{S}=\widehat{R}_{\nu }(S_{W})\widehat{R}_{\nu ^{\prime
}}(S_{W^{\prime }})$ then $\widehat{S}=\widehat{R}_{\nu (S)}(S)$ with 
\begin{equation}
\nu (S)=\nu +\nu ^{\prime }+n-\limfunc{Inert}(M+M^{\prime })  \label{srs}
\end{equation}%
the matrices $M$ and $M^{\prime }$ being associated to $S$ and $S^{\prime }$
by formula (\ref{ms}).
\end{proposition}

\begin{proof}
A straightforward calculation using the composition formula (\ref{tw}) and
the Fresnel integral (\ref{fres}) shows that 
\begin{equation}
\widehat{S}=\left( \frac{1}{2\pi }\right) ^{n}\frac{i^{\nu +\nu ^{\prime }+%
\frac{1}{2}sgn(M+M^{\prime })}}{\sqrt{|\det (S_{W}-I)(S_{W^{\prime
}}-I)(M+M^{\prime })|}}\int e^{\frac{i}{2}\left\langle Nz,z\right\rangle }%
\widehat{T}(z)d^{2n}z  \label{ssm}
\end{equation}%
where $M$ and $M^{\prime }$ correspond to $S_{W}$ and $S_{W^{\prime }}$ by (%
\ref{ms}) and 
\begin{equation*}
N=M-(M+\tfrac{1}{2}J)(M+M^{\prime })^{-1}(M-\tfrac{1}{2}J)\text{.}
\end{equation*}%
We claim that%
\begin{equation}
\det [(S_{W}-I)(S_{W^{\prime }}-I)(M+M^{\prime })]=\det (S-I)  \label{cl1}
\end{equation}%
(hence $M+M^{\prime }$ is indeed invertible), and%
\begin{equation}
N=\tfrac{1}{2}J(S+I)(S-I)^{-1}=M_{S}\text{.}  \label{cl2}
\end{equation}%
Formula (\ref{cl1}) is easy to check by a direct calculation: by definition
of $M$ and $M^{\prime }$ we have, since $\det J=1$,%
\begin{multline*}
\det [(S_{W}-I)(S_{W^{\prime }}-I)(M+M^{\prime })]= \\
\det [(S_{W}-I)(I+(S_{W}-I)^{-1}+(S_{W^{\prime }}-I)^{-1})(S_{W^{\prime
}}-I)]
\end{multline*}%
that is%
\begin{equation*}
\det [(S_{W}-I)(S_{W^{\prime }}-I)(M+M^{\prime })]=\det (S_{W}S_{W^{\prime
}}-I)
\end{equation*}%
which is precisely (\ref{cl1}). Formula (\ref{cl2}) is at first sight more
cumbersome, and one might be tempted to use the oscillator semigroup
calculations of Howe \cite{Howe} at this stage. There is however an easier
way out: assume that $\widehat{S}=\widehat{S}_{W^{\prime \prime },m^{\prime
\prime }}$; we \textit{know} by Proposition \ref{above} that we \textit{must}
have in this case 
\begin{equation*}
N=\tfrac{1}{2}J(S_{W}S_{W^{\prime }}+I)(S_{W}S_{W^{\prime }}-I)^{-1}
\end{equation*}%
and this algebraic identity then holds for all $S=S_{W}S_{W^{\prime }}$
since the free symplectic matrices are dense in $Sp(n)$. Formula (\ref{ssm})
can thus be rewritten%
\begin{equation*}
\widehat{S}=\left( \frac{1}{2\pi }\right) ^{n}\frac{i^{\nu +\nu ^{\prime }+%
\frac{1}{2}sgn(M+M^{\prime })}}{\sqrt{|\det (S-I)|}}\int e^{\frac{i}{2}%
\left\langle M_{S}z,z\right\rangle }\widehat{T}(z)d^{2n}z
\end{equation*}%
and formula (\ref{srs}) follows noting that if $R$ is any real invertible $%
2n\times 2n$ symmetric matrix with $q$ negative eigenvalues we have $\arg
\det R=q\pi $ $\func{mod}2\pi $ and $\frac{1}{2}\limfunc{sgn}R=n-q$ and hence%
\begin{equation*}
\frac{1}{2}sgn(M+M^{\prime })=n-\limfunc{Inert}(M+M^{\prime }).
\end{equation*}
\end{proof}

\section{Concluding Remarks}

We have justified Mehlig and Wilkinson's claim that the metaplectic
operators corresponding to symplectic matrices with no eigenvalues equal to
one can be written in the form (\ref{sf3}); we have in addition shown that
every metaplectic operator can be written as the product of exactly two such
operators. There are however still interesting open problems. It would be
interesting to relate the index\ $\nu $ appearing in (\ref{sf3}) to the
cohomological Maslov index on $Mp(n)$ we constructed in \cite{AIF,Cocycles}:
this would certainly lead to simpler --or at least more tractable--
calculations for the indices intervening in Gutzwiller's formula, as already
demonstrated in our previous paper \cite{MSdG} where we examined the Maslov
index of the monodromy matrix associated to a periodic Hamiltonian orbit. As
shown by (\ref{cz1}) a related mathematical problem would be to express the
Conley--Zehnder index in terms of Leray's index studied in de Gosson \cite%
{JMPA}.

We hope to come back to these important and interesting questions in a near
future, together with applications to various trace formulae.

\begin{acknowledgement}
I wish to thank Keith Hannabuss for having drawn my attention to his work
and having pointed out the relationship between it and the Mehlig--Wilkinson
constructions. I also extend my warmest thanks to the referees for valuable
suggestions and for having pointed out some miscalculations in an early
version of the manuscript.
\end{acknowledgement}

\begin{acknowledgement}
This work has been partially supported by a grant of the Max-Planck-Institut
fuer Gravitationsphysik (Albert-Einstein-Institut, Golm). I wish to thank
Prof. Hermann Nicolai for his kind hospitality.
\end{acknowledgement}

\end{document}